\documentclass{amsart}
\usepackage{amsmath}
\usepackage{amsfonts}

\newtheorem{theorem}{Theorem}
\theoremstyle{plain}

\newtheorem{definition}{Definition}

\newtheorem{lemma}{Lemma}

\newtheorem{remark}{Remark}

\numberwithin{equation}{section}
\input{tcilatex}

\begin{document}
\title{On some nonlinear evolution equation of second order }
\author{Kamal N. Soltanov}
\address{Institute of Math. and Mech. Nat. Acad. of Sci. , Baku, AZERBAIJAN;
Dep. of Math., Fac. of Sci., Hacettepe University, Ankara, TURKEY}
\email{soltanov@hacettepe.edu.tr ; sultan\_kamal@hotmail.com}
\urladdr{}
\date{}

\begin{abstract}
Here we study the abstract nonlinear differential equation of second order
that in special case is the equation of the type of equation of traffic
flow. We prove the solvability theorem for the posed problem under the
appropriate conditions and also investigate the behaviour of the solution.
\end{abstract}

\maketitle

\section{Introduction}

In this article we study the following nonlinear evolution equation 
\begin{equation}
x_{tt}+A\circ F\left( x\right) =g\left( x,A^{-\frac{1}{2}}x_{t}\right) ,\quad t\in \left( 0,T\right) ,\quad 0<T<\infty 
\label{1.1}
\end{equation}%
under the initial conditions 
\begin{equation}
x\left( 0\right) =x_{0},\ x_{t}\left( 0\right) =x_{1}
\label{1.2}
\end{equation}%
here $A$ is a linear operator in a real Hilbert space $H$, $%
F:X\longrightarrow X^{\ast }$ and $g:D\left( g\right) \subseteq H\times
H\longrightarrow H$ are a nonlinear operators, $X$ is a real Banach space.
For example, operator $A$ denotes $-\Delta $ with Dirichlet boundary
conditions (such as homogeneous or periodic) and $f,g$ are functions such as
above, that in the one space dimension case, we can formulate in the form 
\begin{equation}
u_{tt}-\left( f\left( u\right) u_{x}\right) _{x}=g\left( u\right) ,\quad \left( t,x\right) \in R_{+}\times \left( 0,l\right) ,\ l>0,
\label{1.3}
\end{equation}%
\begin{equation}
u\left( 0,x\right) =u_{0}\left( x\right) ,\ u_{t}\left( 0,x\right) =u_{1}\left( x\right) ,\ u\left( t,0\right) =u\left( t,l\right) ,
\label{1.4}
\end{equation}%
where $u_{0}\left( x\right) $, $u_{1}\left( x\right) $ are known functions, $%
f\left( \cdot \right) ,g\left( \cdot \right) :R\longrightarrow R$ are a
continouos functions and $l>0$ is a number. The equation of type (\ref{1.3})
describe mathematical model of the problem from theory of the flow in
networks as is affirmed in articles \cite{a-r}, \cite{b-c-g-h-p}, \cite{c-t}%
, \cite{h-r} (e. g. Aw-Rascle equations, Antman--Cosserat model, etc.). As
in the survey \cite{b-c-g-h-p} is noted such a study can find application in
accelerating missiles and space crafts, components of high-speed machinery,
manipulator arm, microelectronic mechanical structures, components of
bridges and other structural elements. Balance laws are hyperbolic partial
differential equations that are commonly used to express the fundamental
dynamics of open conservative systems (e.g. \cite{c-t}). As the survey \cite%
{b-c-g-h-p} possess of the sufficiently exact explanations of the
significance of equations of such type therefore we not stop on this theme.

This article is organized as follows. In the section 2 we study the
solvability of the nonlinear equation of second order in the Banach spaces,
for which we found the sufficient conditions and proved the existence
theorem. In the section 3 we investigate the global behaviour of solutions
of the posed problem.

\section{Solvability of problem (\protect\ref{1.1}) - (\protect\ref{1.2})}

Let $A$ is a symmetric linear operator densely defined in a real Hilbert
space $H$ and positive, $A$ has a self-adjoint extension, moreover there is
linear operator $B$ defined in $H$ such that $A\equiv B^{\ast }\circ B$,
here $f:R\longrightarrow R$ is continuous as function, $X$ is a real
reflexive Banach space and $X\subset H$, $g:D\left( g\right) \subseteq
H\times H\longrightarrow H$, where $g:R^{2}\longrightarrow R$ is a
continuous as function and $x:\left[ 0,T\right) \longrightarrow X$ is an
unknown function. Let $F\left( r\right) $ as a function is defined as $%
F\left( r\right) =\underset{0}{\overset{r}{\int }}f\left( s\right) ds$. Let
the inequation $\left\Vert x\right\Vert _{H}\leq \left\Vert Bx\right\Vert
_{H}$ is valid for any $x\in D\left( B\right) $. We denote by $V$, $W$ and
by $Y$ the spaces defined as $V\equiv \left\{ y\in H\left\vert \ By\in
H\right. \right\} $, $W=\left\{ x\in H\left\vert \ Ax\in H\right. \right\} $
and as $Y\equiv \left\{ x\in X\left\vert \ Ax\in X\right. \right\} $,
respectively, for which inclusions $W\subset V\subset H$ are compact and $%
Y\subset W$.

Let $H$ is the real separable Hilbert space, $X$ is the reflexive Banach
space and $X\subset H\subset X^{\ast }$; $V$ is the previously defined
space. It is clear that $W\subset V\subset H\subset V^{\ast }\subset W^{\ast
}$ are a framed spaces by $H$, these inclusions are compact and $X\subset
V^{\ast }$ then one can define the framed spaces $Y\subset V\subset H\subset
V^{\ast }\subset Y^{\ast }$ then $X\subset V^{\ast }\subset Y^{\ast }$ are
compact, with use the property of the operator $A$. Assume that operator $A$
such as $A:$ $V_{B}\longrightarrow V_{B}^{\ast }$ and $A:X^{\ast
}\longrightarrow Y^{\ast }$. Consequently, we get $A\circ F:X\longrightarrow
Y^{\ast }$ and $A\circ F\circ A:Y\longrightarrow Y^{\ast }$. Moreover we
assume that $\left[ X^{\ast },Y\right] _{\frac{1}{2}}\subseteq V$.

Since operator $A$ is invertible, here one can set the function $y\left(
t\right) =A^{-1}x\left( t\right) $ for any $t\in \left( 0,T\right) $, in the
other words one can assume the denotation $x\left( t\right) =Ay\left(
t\right) $.

We will interpret the solution of the problem (\ref{1.1}) - (\ref{1.2}) by
the following manner.

\begin{definition}
\label{D-1}A function $x:$ $\left( 0,T\right) \longrightarrow X$, $x\in
C^{0}\left( 0,T;X\right) \cap C^{1}\left( 0,T;V^{\ast }\right) \cap
C^{2}\left( 0,T;Y^{\ast }\right) $, $x=Ay$, is called a weak solution of
problem (\ref{1.1}) - (\ref{1.2}) if $x$ a. e. $t\in \left( 0,T\right) $
satisfies the following equation 
\begin{equation}
\frac{d^{2}}{dt^{2}}\left\langle x,z\right\rangle +\left\langle A\circ
F\left( x\right) ,z\right\rangle =\left\langle g\left( x,By_{t}\right)
,z\right\rangle  \label{2.1}
\end{equation}%
for any $z\in Y$ and the initial conditions (\ref{1.2}) (here and farther
the expression $\left\langle \cdot ,\cdot \right\rangle $ denotes the dual
form for the pair: the Banach space and his dual).
\end{definition}

Consider the following conditions

(\textit{i}) Let $A:$ $W\subset H\longrightarrow H$ is the selfadjoint and
positive operator, moreover $A:$ $V\longrightarrow V^{\ast }$, $A:X^{\ast
}\longrightarrow Y^{\ast }$, there exists an linear operator $%
B:V\longrightarrow H$ that satisfies the equation $Ax\equiv \left( B^{\ast
}\circ B\right) x$ for any $x\in D\left( A\right) $ and $\left\Vert
x\right\Vert _{H}\leq \left\Vert Bx\right\Vert _{H}=\left\Vert x\right\Vert
_{V}$.

(\textit{ii}) Let $F:X\longrightarrow X^{\ast }$ is the continuously
differentiable and monotone operator with the potential $\Phi $ that is the
functional defined on $X$ (his Frechet derivative is the operator $F$).
Moreover for any $x\in X$ the following inequalities hold 
\begin{equation*}
\left\Vert F\left( x\right) \right\Vert _{X^{\ast }}\leq a_{0}\left\Vert
x\right\Vert _{X}^{p-1}+a_{1}\left\Vert x\right\Vert _{H};\quad \left\langle
F\left( x\right) ,x\right\rangle \geq b_{0}\left\Vert x\right\Vert
_{X}^{p}+b_{1}\left\Vert x\right\Vert _{H}^{2},
\end{equation*}%
where $a_{0},b_{0}>0$, $a_{1},b_{1}\geq 0$, $p>2$ are numbers.

(\textit{iii}) Assume $g:H\times V\longrightarrow H$ is a continuous
operator that satisfies the condition 
\begin{equation*}
\left\vert \left\langle g\left( x,y\right) -g\left( x_{1},y_{1}\right)
,z\right\rangle \right\vert \leq g_{1}\left\vert \left\langle
x-x_{1},z\right\rangle \right\vert +g_{2}\left\vert \left\langle
y-y_{1},z\right\rangle \right\vert ,
\end{equation*}%
for any $\left( x,y\right) ,\left( x_{1},y_{1}\right) \in H\times H$, $z\in
H $ and consequently for any $\left( x,y\right) \in H\times H$ the
inequation 
\begin{equation*}
\left\Vert g\left( x,y\right) \right\Vert _{H}\leq g_{1}\left\Vert
x\right\Vert _{H}+g_{2}\left\Vert y\right\Vert _{H}+g_{0},\quad g_{0}\geq
\left\Vert g\left( 0,0\right) \right\Vert _{H}
\end{equation*}%
holds, where $g_{0}$ is a number.

In the beginning for the investigation of the posed problem we set the
following expression in order to obtain of the a priori estimations 
\begin{equation*}
\left\langle x_{tt},y_{t}\right\rangle +\left\langle A\circ F\left( x\right)
,y_{t}\right\rangle =\left\langle g\left( x,By_{t}\right) ,y_{t}\right\rangle
\end{equation*}%
here element $y$ is defined as the solution of the equation $Ay\left(
t\right) =x\left( t\right) $, i.e. $y\left( t\right) =A^{-1}x\left( t\right) 
$ for any $t\in \left( 0,T\right) $ as was already mentioned above.

Hence follow 
\begin{equation*}
\left\langle By_{tt},By_{t}\right\rangle +\left\langle F\left( x\right)
,x_{t}\right\rangle =\left\langle g\left( x,By_{t}\right)
,y_{t}\right\rangle ,
\end{equation*}%
or 
\begin{equation}
\frac{1}{2}\frac{d}{dt}\left\Vert By_{t}\right\Vert _{H}^{2}+\frac{d}{dt}\Phi \left( x\right) =\left\langle g\left( x,By_{t}\right) ,y_{t}\right\rangle ,
\label{2.2}
\end{equation}%
where $\Phi \left( x\right) $ is the functional that defined as $\Phi \left(
x\right) =\underset{0}{\overset{1}{\int }}\left\langle F\left( sx\right)
,x\right\rangle ds$ (see, \cite{le}).

Then using condition (\textit{iii}) on $g\left( x,By_{t}\right) $ in (\ref%
{2.2}) one can obtain 
\begin{equation*}
\frac{1}{2}\frac{d}{dt}\left\Vert By_{t}\right\Vert _{H}^{2}+\frac{d}{dt}%
\Phi \left( x\right) \leq \left\Vert g\left( x,By_{t}\right) \right\Vert
_{H}^{2}+\left\Vert y_{t}\right\Vert _{H}^{2}\leq
\end{equation*}%
\begin{equation*}
2\left( g_{1}^{2}\left\Vert x\right\Vert _{H}^{2}+g_{2}^{2}\left\Vert
By_{t}\right\Vert _{H}^{2}+g_{0}^{2}\right) +\left\Vert y_{t}\right\Vert
_{H}^{2}\leq \widetilde{C}\left( \left\Vert x\right\Vert _{H}^{2}+\frac{1}{2}%
\left\Vert By_{t}\right\Vert _{H}^{2}+g_{0}^{2}\right)
\end{equation*}%
\ here one can use the estimation $\left\Vert x\right\Vert _{H}^{2}\leq 
\widetilde{c}\left( \Phi \left( x\right) +1\right) $ (if $b_{1}>0$ then $%
\left\Vert x\right\Vert _{H}^{2}\leq \widetilde{c}\Phi \left( x\right) $) as 
$2<p$ by virtue of the condition \textit{(ii)}. Consequently we get to the
Cauchy problem for the inequation 
\begin{equation}
\frac{d}{dt}\left( \frac{1}{2}\left\Vert By_{t}\right\Vert _{H}^{2}\left( t\right) +\Phi \left( x\left( t\right) \right) \right) \leq C_{0}\left( \frac{1}{2}\left\Vert By_{t}\right\Vert _{H}^{2}\left( t\right) +\Phi \left( x\left( t\right) \right) \right) +C_{1}
\label{2.3}
\end{equation}%
with the initial conditions 
\begin{equation}
x\left( t\right) \left\vert \ _{t=0}\right. =x_{0};\quad y_{t}\left( t\right) \left\vert \ _{t=0}\right. =A^{-1}x_{t}\left\vert \ _{t=0}\right. =A^{-1}x_{1}
\label{2.4}
\end{equation}%
where $C_{j}\geq 0$ are constants independent of $x\left( t\right) $. \ From
here follows 
\begin{equation*}
\frac{1}{2}\left\Vert By_{t}\right\Vert _{H}^{2}\left( t\right) +\Phi \left(
x\left( t\right) \right) \leq e^{tC_{0}}\left[ \left\Vert By_{1}\right\Vert
_{2}^{2}+2\Phi \left( x_{0}\right) \right] +\frac{C_{1}}{C_{0}}\left(
e^{tC_{0}}-1\right) .
\end{equation*}

This give to we the following estimations for every $T\in \left( 0,\infty
\right) $ 
\begin{equation}
\left\Vert By_{t}\right\Vert _{H}^{2}\left( t\right) \leq C\left( x_{0},x_{1}\right) e^{C_{0}T},\quad \Phi \left( x\left( t\right) \right) \leq C\left( x_{0},x_{1}\right) e^{C_{0}T},
\label{2.5}
\end{equation}%
for a. e. $t\in \left( 0,T\right) $, i.e. $y=A^{-1}x$ is contained in the
bounded subset of the space $y\in C^{1}\left( 0,T;V\right) \cap C^{0}\left(
0,T;Y\right) $, consequently we obtain that if the weak solution $x\left(
t\right) $ exists then it belong to a bounded subset of the space $%
C^{0}\left( 0,T;X\right) \cap C^{1}\left( 0,T;V^{\ast }\right) $.

Hence one can wait, that the following inclusion 
\begin{equation*}
y\in C^{2}\left( 0,T;X^{\ast }\cap H\right) \cap C^{1}\left( 0,T;X^{\ast
}\cap V\right) \cap C^{0}\left( 0,T;Y\right)
\end{equation*}%
holds by virtue of (\ref{2.1}) in the assumption that $x=Ay$ is a solution
of the posed problem in the sense of Definition \ref{D-1}.

In order to prove of the solvability theorem we will use the Galerkin
approach. Let the system $\left\{ y^{k}\right\} _{k=1}^{\infty }\subset Y$
be total in $Y$ such that it is complete in the spaces $Y,V$, and also in
the spaces $X,H$. We will seek out of the approximative solutions $%
y_{m}\left( t\right) $, consequently and $x_{m}\left( t\right) $, in the
form 
\begin{equation*}
x_{m}\left( t\right) \equiv Ay_{m}\left( t\right) =\overset{m}{\underset{k=1}%
{\sum }}c_{i}\left( t\right) Ay^{k}\text{ or }x_{m}\left( t\right) \in
span\left\{ y^{1},...,y^{m}\right\}
\end{equation*}%
as the solutions of the problem locally with respect to $t$, where $%
c_{i}\left( t\right) $ are as the unknown functions that will be defined as
solutions of the following Cauchy problem for system of ODE 
\begin{equation*}
\frac{d^{2}}{dt^{2}}\left\langle x_{m},y^{j}\right\rangle +\left\langle
F\left( x_{m}\right) ,Ay^{j}\right\rangle =\left\langle g\left(
x_{m},By_{mt}\right) ,y^{j}\right\rangle ,\quad j=1,2,...,m
\end{equation*}%
\begin{equation*}
x_{m}\left( 0\right) =x_{0m},\ x_{tm}\left( 0\right) =x_{1m},
\end{equation*}%
where $x_{0m}$ and $x_{1m}$ are contained in $span\left\{
y^{1},...,y^{m}\right\} $, $m=1,2,...$, moreover 
\begin{equation*}
x_{0m}\longrightarrow x_{0}\quad \text{in \ }\left[ X,Y\right] _{\frac{1}{2}%
}\subseteq V;\quad x_{1m}\longrightarrow x_{1}\quad \text{in \ }X,\
m\nearrow \infty .
\end{equation*}

Thus we obtain the following problem 
\begin{equation}
\frac{d^{2}}{dt^{2}}\left\langle x_{m},y^{j}\right\rangle +\left\langle
F\left( x_{m}\right) ,Ay^{j}\right\rangle =\left\langle g\left(
x_{m},By_{mt}\right) ,y^{j}\right\rangle ,\quad j=1,2,...,m  \label{2.6}
\end{equation}%
\begin{equation*}
\left\langle x_{m}\left( t\right) ,y^{j}\right\rangle \left\vert \
_{t=0}\right. =\left\langle x_{0m},y^{j}\right\rangle ,\ \frac{d}{dt}%
\left\langle x_{m}\left( t\right) ,y^{j}\right\rangle \left\vert \
_{t=0}\right. =\left\langle x_{1m},y^{j}\right\rangle
\end{equation*}%
that solvable by virtue of estimates (\ref{2.5}) on $\left( 0,T\right) $ for
any $m=1,2,...$, $j=1,2,...$ and $T>0.$ Hence we set 
\begin{equation}
\frac{d^{2}}{dt^{2}}\left\langle x_{m},z\right\rangle +\left\langle F\left(
x_{m}\right) ,Az\right\rangle =\left\langle g\left( x_{m},By_{mt}\right)
,z\right\rangle  \label{2.7}
\end{equation}%
for any $z\in Y$ and $m=1,2,...$.

Consequently with use of the known procedure (\cite{li}, \cite{s1}, \cite{z}%
) we obtain, $y_{mt}\in C^{0}\left( 0,T;V\right) $, $y_{m}\in C^{0}\left(
0,T;Y\right) $ and $x_{m}\in C^{0}\left( 0,T;X\right) $, $\ x_{mt}\in
C^{0}\left( 0,T;V^{\ast }\right) $, moreover they are contained in the
bounded subset of these spaces for any $m=1,2,...$. Hence from (\ref{2.5})
we get 
\begin{equation*}
x_{mtt}\in C^{0}\left( 0,T;Y^{\ast }\right) \text{ \ \ or \ \ }x_{m}\in
C^{2}\left( 0,T;Y^{\ast }\right) \text{, (}V^{\ast }\subset Y^{\ast }\text{).%
}
\end{equation*}%
Thus we obtain, that the sequence $\left\{ x_{m}\right\} _{m=1}^{\infty }$
of the approximated solutions of the problem is contained in a bounded
subset of the space 
\begin{equation*}
C^{0}\left( 0,T;X\right) \cap C^{1}\left( 0,T;V^{\ast }\right) \cap
C^{2}\left( 0,T;Y^{\ast }\right)
\end{equation*}%
or $\left\{ x_{m}\right\} _{m=1}^{\infty }$ such that for a. e. $t\in \left(
0,T\right) $ takes place the following inclusions $\left\{ y_{m}\left(
t\right) \right\} _{m=1}^{\infty }\subset Y\subset X\subset H$, $\left\{
y_{mt}\left( t\right) \right\} _{m=1}^{\infty }\subset V$, $\left\{
y_{mtt}\left( t\right) \right\} _{m=1}^{\infty }\subset X^{\ast }$. So we
have 
\begin{equation*}
\left\{ y_{m}\left( t\right) \right\} _{m=1}^{\infty }\subset C^{0}\left(
0,T;Y\right) \cap C^{1}\left( 0,T;V\right) \cap C^{2}\left( 0,T;X^{\ast
}\right)
\end{equation*}%
therefore $\left\{ y_{m}\left( t\right) \right\} _{m=1}^{\infty }$ possess a
precompact subsequence in $C^{1}\left( 0,T;\left[ X^{\ast },Y\right] _{\frac{%
1}{2}}\right) $ and in $C^{1}\left( 0,T;V\right) $, as $\left[ X^{\ast },Y%
\right] _{\frac{1}{2}}\subseteq V$ by virtue of conditions on $X$ and $A$
(by virtue of well known results, see, e. g. \cite{b-s}, \cite{s} etc.).
From here follows $y_{m}\left( t\right) \longrightarrow y\left( t\right) $
in $C^{1}\left( 0,T;V\right) $ for $m\nearrow \infty $ (Here and hereafter
in order to abate the number of index we don't changing of indexes of
subsequences). Then the sequence $\left\{ F\left( Ay_{m}\left( t\right)
\right) \right\} _{m=1}^{\infty }\subset X^{\ast }$ and bounded for a. e. $%
t\in \left( 0,T\right) $; the sequence 
\begin{equation*}
\left\{ g\left( x_{m}\left( t\right) ,x_{mt}\left( t\right) \right) \right\}
_{m=1}^{\infty }\equiv \left\{ g\left( Ay_{m}\left( t\right) ,By_{mt}\left(
t\right) \right) \right\} _{m=1}^{\infty }\subset H
\end{equation*}%
and bounded for a. e. $t\in \left( 0,T\right) $ also, by virtue of the
condition \textit{(iii).} Indeed, for any $m$ the estimation 
\begin{equation*}
\left\Vert g\left( Ay_{m},By_{mt}\right) \right\Vert _{H}\left( t\right)
\leq \left\Vert Ay_{m}\left( t\right) \right\Vert _{H}+\left\Vert
By_{mt}\left( t\right) \right\Vert _{H}+\left\Vert g\left( 0,0\right)
\right\Vert _{H}
\end{equation*}%
holds and therefore $\left\{ g\left( Ay_{m}\left( t\right) ,By_{mt}\left(
t\right) \right) \right\} _{m=1}^{\infty }$ is contained in a bounded subset
of $H$ for a. e. $t\in \left( 0,T\right) $. Consequently $\left\{ F\left(
Ay_{m}\right) \right\} _{m=1}^{\infty }$ and $\left\{ g\left( Ay_{m}\left(
t\right) ,By_{mt}\left( t\right) \right) \right\} _{m=1}^{\infty }$ have an
weakly converging subsequences to $\eta \left( t\right) $ and $\theta \left(
t\right) $ in $X^{\ast }$ and $H$, respectively, for a. e. $t\in \left(
0,T\right) $. Hence one can pass to the limit in (\ref{2.7}) with respect to 
$m\nearrow \infty $. Then we obtain the following equation 
\begin{equation}
\frac{d^{2}}{dt^{2}}\left\langle x,z\right\rangle +\left\langle A\eta \left(
t\right) ,z\right\rangle =\left\langle \theta \left( t\right)
,z\right\rangle .  \label{2.8}
\end{equation}

So for us is remained to show the following: if the sequence $\left\{
x_{m}\left( t\right) \right\} _{m=1}^{\infty }\equiv \left\{ Ay_{m}\left(
t\right) \right\} _{m=1}^{\infty }$ is weakly converging to $x\left(
t\right) =Ay\left( t\right) $ then $\eta \left( t\right) =F\left( x\left(
t\right) \right) $ and $\theta \left( t\right) =g\left( x\left( t\right)
,By_{t}\left( t\right) \right) $. In order to show these equations are
fulfilled we will use the monotonicity of $F$ and the condition \textit{%
(iii). }

We start to show $\theta \left( t\right) =g\left( \left( t\right)
,By_{t}\left( t\right) \right) $ as $x\in X\subset H$ and $y_{t}\in V$, $%
By_{t}\in H$ therefore $g\left( x,By_{t}\right) $ is defined for a. e. $t\in
\left( 0,T\right) $. Consequently one can consider of the expression 
\begin{equation*}
\left\langle g\left( Ay_{m}\left( t\right) ,By_{mt}\left( t\right) \right)
-g\left( Ay\left( t\right) ,By_{t}\left( t\right) \right) ,\widehat{y}%
\right\rangle
\end{equation*}%
for any $\widehat{y}\in C^{0}\left( 0,T;Y\right) \cap C^{1}\left(
0,T;V\right) $. So we set this expression and investigate this for any $%
\widehat{y}\in C^{0}\left( 0,T;Y\right) \cap C^{1}\left( 0,T;V\right) $ then
we have 
\begin{equation*}
\left\vert \left\langle g\left( Ay_{m}\left( t\right) ,By_{mt}\left(
t\right) \right) -g\left( Ay\left( t\right) ,By_{t}\left( t\right) \right) ,%
\widehat{y}\right\rangle \right\vert \leq
\end{equation*}%
\begin{equation}
g_{1}\left\vert \left\langle Ay_{m}\left( t\right) -Ay\left( t\right) ,%
\widehat{y}\left( t\right) \right\rangle \right\vert +g_{2}\left\vert
\left\langle By_{mt}\left( t\right) -By_{t}\left( t\right) ,\widehat{y}%
\left( t\right) \right\rangle \right\vert  \label{2.9}
\end{equation}%
that takes place by virtue of the condition \textit{(iii).} Using here the
weak convergences of $Ay_{m}\left( t\right) \rightharpoonup Ay\left(
t\right) $ and $By_{mt}\left( t\right) \longrightarrow By_{t}\left( t\right) 
$ and by passing to the limit in the inequatin (\ref{2.9}) with respect to $%
m:$ $m\nearrow \infty $ we get 
\begin{equation*}
\left\vert \left\langle \theta \left( t\right) -g\left( Ay\left( t\right)
,By_{t}\left( t\right) \right) ,\widehat{y}\right\rangle \right\vert \leq 0
\end{equation*}%
for any $\widehat{y}\in C^{0}\left( 0,T;H\right) $. Consequently the
equation $\theta \left( t\right) =g\left( \left( t\right) ,By_{t}\left(
t\right) \right) $ holds, then the following equation is valid 
\begin{equation*}
\frac{d^{2}}{dt^{2}}\left\langle x,z\right\rangle +\left\langle A\eta \left(
t\right) ,z\right\rangle =\left\langle g\left( x\left( t\right)
,By_{t}\left( t\right) \right) ,z\right\rangle
\end{equation*}%
for any $z\in Y$, \ as $\left\{ y^{k}\right\} _{k=1}^{\infty }$ is complete
in $Y$\ that display fulfilling of equation 
\begin{equation}
A\eta \left( t\right) =g\left( x\left( t\right) ,By_{t}\left( t\right)
\right) -\frac{d^{2}x}{dt^{2}}  \label{2.10}
\end{equation}
in the sense of $Y^{\ast }$.

In order to show the equation $\eta \left( t\right) =F\left( x\left(
t\right) \right) $ one can use the monotonicity of $F$. So the following
inequation holds 
\begin{equation*}
\left\langle A\circ F\left( Az\right) -A\circ F\left( Ay\right)
,z-y\right\rangle =\left\langle F\left( Az\right) -F\left( Ay\right)
,Az-Ay\right\rangle =
\end{equation*}%
\begin{equation*}
\left\langle F\left( \widetilde{x}\right) -F\left( x\right) ,\widetilde{x}%
-x\right\rangle \geq 0
\end{equation*}%
for any $y,z\in Y$, $Ay=x$ and $Az=\widetilde{x}$ by condition \textit{(i)}.
Then one can write 
\begin{equation*}
0\leq \left\langle F\left( x_{m}\right) -F\left( \widetilde{x}\right) ,x_{m}-%
\widetilde{x}\right\rangle =\left\langle F\left( Ay_{m}\right) -F\left(
Az\right) ,Ay_{m}-Az\right\rangle =
\end{equation*}%
take account here the equation (\ref{2.5}) \ 
\begin{equation*}
\left\langle F\left( Ay_{m}\right) ,Ay_{m}\right\rangle -\left\langle \frac{%
d^{2}}{dt^{2}}x_{m}-g\left( x_{m},By_{mt}\right) ,z\right\rangle
-\left\langle F\left( Az\right) ,Ay_{m}-Az\right\rangle =
\end{equation*}%
\begin{equation}
\left\langle F\left( x_{m}\right) ,x_{m}\right\rangle -\left\langle \frac{%
d^{2}}{dt^{2}}x_{m}-g\left( x_{m},By_{mt}\right) ,z\right\rangle
-\left\langle F\left( \widetilde{x}\right) ,x_{m}-\widetilde{x}\right\rangle
.  \label{2.11}
\end{equation}%
Here one can use the well-known inequation 
\begin{equation*}
\lim \sup \left\langle F\left( x_{m}\right) ,x_{m}\right\rangle \leq
\left\langle \eta ,x\right\rangle =\left\langle \eta ,Ay\right\rangle
=\left\langle A\eta ,y\right\rangle .
\end{equation*}%
Then passing to the limit in (\ref{2.11}) with respect to $m:m\nearrow
\infty $ we obtain 
\begin{equation*}
0\leq \left\langle A\eta ,y\right\rangle -\left\langle \frac{d^{2}}{dt^{2}}%
x-g\left( x,By_{t}\right) ,z\right\rangle -\left\langle F\left( \widetilde{x}%
\right) ,x-\widetilde{x}\right\rangle =
\end{equation*}%
\begin{equation*}
\left\langle A\eta ,y\right\rangle -\left\langle A\eta ,z\right\rangle
-\left\langle F\left( \widetilde{x}\right) ,Ay-Az\right\rangle =\left\langle
A\eta -A\circ F\left( \widetilde{x}\right) ,y-z\right\rangle
\end{equation*}%
by virtue of (\ref{2.10})

Consequently we obtain, that the equation $A\eta \left( t\right) =A\circ
F\left( x\right) $ holds since $z$ is arbitrary element of $Y$.

Now for us is remained to show the obtained function $x(t)=Ay(t)$ satisfy
the initial conditions. Consider the following equation 
\begin{equation*}
\left\langle y_{mt},Ay_{m}\right\rangle \left( t\right) =\underset{0}{%
\overset{t}{\int }}\left\langle \frac{d^{2}}{ds^{2}}Ay_{m},y_{m}\right%
\rangle ds+\underset{0}{\overset{t}{\int }}\left\langle \frac{d}{ds}By_{m},%
\frac{d}{ds}By_{m}\right\rangle ds+\left\langle y_{1m},Ay_{0m}\right\rangle =
\end{equation*}%
\begin{equation*}
\underset{0}{\overset{t}{\int }}\left\langle \frac{d^{2}}{ds^{2}}%
y_{m},Ay_{m}\right\rangle ds+\underset{0}{\overset{t}{\int }}\left\Vert 
\frac{d}{ds}By_{m}\right\Vert _{H}^{2}ds+\left\langle
y_{1m},Ay_{0m}\right\rangle
\end{equation*}%
for $m=1,2,...$, here $x_{m}(t)=Ay_{m}(t)$. Hence we get: the left side is
bounded as far as all addings items in the right side are bounded by virtue
of the obtained estimations. Therefore one can pass to limit with respect to 
$m$ as here $y_{mt}$ are continous with respect to $t$ for any $m$ then $%
y_{mt}$ strongly converges to $y_{t}$ and $Ay_{m}$ weakly converges to $Ay$
in $H$. It must be noted the equation 
\begin{equation*}
\underset{m\longrightarrow \infty }{\lim }\overset{t}{\underset{0}{\int }}%
\left\Vert \frac{d}{ds}By_{m}\right\Vert _{H}^{2}dxds=\underset{0}{\overset{t%
}{\int }}\left\Vert \frac{d}{ds}By\right\Vert _{H}^{2}ds
\end{equation*}%
holds by virtue of the above reasonings that $\left\{ y_{m}\left( t\right)
\right\} _{m=1}^{\infty }$ is a precompact subset in $C^{1}\left(
0,T;V\right) $. Consequently the left side converges to the expression of
such type, i.e. to $\left\langle y_{t},Ay\right\rangle \left( t\right) $.

The obtained results shows that the following convergences are just: $%
x_{m}\left( t\right) =Ay_{m}\left( t\right) \rightharpoonup Ay\left(
t\right) =x\left( t\right) $ in $X$, $x_{mt}\left( t\right)
=Ay_{mt}\rightharpoonup Ay_{t}=x_{t}\left( t\right) $ in $V^{\ast }$. From
here follows, that the initial conditions are fulfilled in the sense of $X$
and $V^{\ast }$, respectively.

Thus the following existence resut is proved.

\begin{theorem}
Let spaces $H,V,W,X,Y$ that are defined above satisfy all above mentioned
conditions and condition (i)-(iii) are fulfilled then problem (\ref{1.1}) - (%
\ref{1.2}) is solvable in the space $C^{0}\left( 0,T;X\right) \cap
C^{1}\left( 0,T;V\right) \cap C^{2}\left( 0,T;Y^{\ast }\right) $ for any $%
x_{0}\in V\cap \left[ X^{\ast },Y\right] _{\frac{1}{2}}$ and $x_{1}\in H$ in
the sense of Definition \ref{D-1}.
\end{theorem}

\begin{remark}
This theorem shows that there exist a flow $S\left( t\right) $ defined in $%
V\times X$ and the solution of the problem (\ref{1.1}) - (\ref{1.2}) one can
represent as $x\left( t\right) =S\left( t\right) \circ \left(
x_{0},x_{1}\right) $.
\end{remark}

\section{Behaviour of solutions of problem (\protect\ref{1.1}) - (\protect
\ref{1.2})}

Here we consider problem under the following complementary conditions:

(\textit{iv}) Let $g\left( x,By_{t}\right) =0$ and $\left\Vert x\right\Vert
_{H}^{p}\left( t\right) \leq c_{0}\Phi \left( x\left( t\right) \right) $ for
some $c_{0}>0$.

We set a function $E(t)=\Vert Bw\Vert _{H}^{2}(t)$ and consider this
function on the solution of problem (\ref{1.1}) - (\ref{1.2}), then\ for $%
E(t)=\Vert By\Vert _{H}^{2}(t)$ we have 
\begin{equation}
E%
{\acute{}}%
\left( t\right) =2\left\langle By_{t},By\right\rangle \leq \left\Vert
By_{t}\right\Vert _{2}^{2}\left( t\right) +\left\Vert By\right\Vert
_{2}^{2}\left( t\right) ,  \label{3.1}
\end{equation}%
where $y=A^{-1}x$. Here we will use equation (\ref{2.5}). For this we lead
the following equation 
\begin{equation*}
\frac{1}{2}\left\Vert By_{s}\right\Vert _{H}^{2}\left( s\right) +\Phi \left(
x\left( s\right) \right) \left\vert _{0}^{t}\ \right. =0
\end{equation*}%
as $g\left( x,By_{t}\right) =0$.\footnote{%
We would like to note that this equation shows the stability of the energy
of the considered system in this case.} Hence 
\begin{equation*}
\frac{1}{2}\left\Vert By_{s}\right\Vert _{H}^{2}\left( t\right) +\Phi \left(
x\left( t\right) \right) =\frac{1}{2}\left\Vert By_{1}\right\Vert
_{H}^{2}\left( t\right) +\Phi \left( x_{0}\right)
\end{equation*}%
and 
\begin{equation*}
\left\Vert By_{t}\right\Vert _{H}^{2}\left( t\right) =-2\Phi \left( x\left(
t\right) \right) +\left\Vert By_{1}\right\Vert _{H}^{2}+2\Phi \left(
x_{0}\right) .
\end{equation*}

Granting this in (\ref{3.1}) we get 
\begin{equation*}
E%
{\acute{}}%
\left( t\right) \leq E\left( t\right) -E^{r}\left( t\right) +\left\Vert
By_{1}\right\Vert _{H}^{2}+2\Phi \left( x_{0}\right)
\end{equation*}%
by virtue of the condition $\Phi \left( x\right) \geq c_{0}\left\Vert
x\right\Vert _{X}^{p}$ and of the continuity of embedding $X\subset H$, $%
r=p/2$.

So denoted by $z\left( t\right) =E\left( t\right) $ we have the Cauchy
problem for differential inequality 
\begin{equation}
z%
{\acute{}}%
\left( t\right) \leq z\left( t\right) -cz^{r}\left( t\right) +C\left(
x_{0},x_{1}\right) ,\quad z\left( 0\right) =\left\Vert By_{0}\right\Vert
_{H}^{2},  \label{3.2}
\end{equation}%
that we will investigate. Inequation (\ref{3.2}) one can rewrite in the form 
\begin{equation*}
\left( z\left( t\right) +kC\left( x_{0},x_{1}\right) \right) 
{\acute{}}%
\leq z\left( t\right) +kC\left( x_{0},x_{1}\right) -\delta \left[ z\left(
t\right) +kC\left( x_{0},x_{1}\right) \right] ^{r},
\end{equation*}%
where $k>1$ is a number and $\delta =\delta \left( c,C,k,r\right) >0$ is
sufficiently small number. Then solving this problem we get 
\begin{equation*}
z\left( t\right) +kC\left( x_{0},x_{1}\right) \leq \left[ e^{\left(
1-r\right) t}\left( z_{0}+kC\left( x_{0},x_{1}\right) \right) ^{1-r}+\delta
\left( 1-e^{\left( 1-r\right) t}\right) \right] ^{\frac{1}{1-r}}
\end{equation*}%
or 
\begin{equation*}
E\left( t\right) \leq \left[ e^{\left( 1-r\right) t}\left( \left\Vert
By_{0}\right\Vert _{H}^{2}+kC\left( x_{0},x_{1}\right) \right) ^{1-r}+\delta
\left( 1-e^{\left( 1-r\right) t}\right) \right] ^{\frac{1}{1-r}}-kC\left(
x_{0},x_{1}\right)
\end{equation*}%
\begin{equation}
\Vert By\Vert _{H}^{2}(t)\leq \frac{e^{t}\left( \left\Vert By_{0}\right\Vert
_{H}^{2}+kC\left( x_{0},x_{1}\right) \right) }{\left[ 1+\delta \left(
\left\Vert By_{0}\right\Vert _{H}^{2}+kC\left( x_{0},x_{1}\right) \right)
^{r-1}\left( e^{\left( r-1\right) t}-1\right) \right] ^{\frac{1}{r-1}}}%
-kC\left( x_{0},x_{1}\right) .  \tag{3.3}
\end{equation}%
here the right side is greater than zero, because $\delta \leq \frac{k-1}{%
k^{r}C^{r}}$ and $2r=p>2$.

Thus is proved the result

\begin{lemma}
Under conditions (i), (ii), (iv) the function $y(t)$, defined by the
solution of problem (\ref{2.1})-(\ref{2.2}), for any $t>0$ is contained in
ball $B_{l}^{X\cap V}\left( 0\right) \subset X\cap V$ depending from the
initial values $\left( x_{0},x_{1}\right) \in \left( X\cap V\right) \times H$%
, here $l=l\left( x_{0},x_{1},p\right) >0$. \ 
\end{lemma}

\end{document}